\documentclass[review]{elsarticle}

\usepackage{lineno,hyperref}
\usepackage{microtype}
\usepackage{amsmath}
\usepackage{amssymb,amsfonts,amsthm,latexsym,mathrsfs}
\usepackage{eucal}
\usepackage{color}

\theoremstyle{plain}

\newtheorem*{theorem*}{{\bf Theorem}}

\newtheorem*{corollary*}{{\bf Corollary}}

\theoremstyle{definition}

\theoremstyle{remark}

\numberwithin{equation}{subsection}
\modulolinenumbers[5]

\journal{Journal of Pure and Applied Algebra}

\begin{document}

\begin{frontmatter}

\title{An exact upper bound for sums of element orders in non-cyclic finite groups}

\author[1]{Marcel Herzog}
\ead{herzogm@post.tau.ac.il}

\author[2]{Patrizia Longobardi\corref{mycorrespondingauthor}}
\cortext[mycorrespondingauthor]{corresponding author}
\ead{plongobardi@unisa.it}

\author[2]{Mercede Maj}
\ead{mmaj@unisa.it}

\address[1]{Department of Mathematics,
Raymond and Beverly Sackler Faculty of Exact Sciences,

Tel Aviv University,  Tel Aviv, Israel}
\address[2]{Dipartimento di Matematica,
Universit\`a di Salerno,
via Giovanni Paolo II, 132, 84084 Fisciano (Salerno), Italy}

\begin{abstract}  
Denote the sum of element orders in a finite group $G$
by $\psi(G)$ and let $C_n$ denote the cyclic group of order $n$.
Suppose that $G$ is a non-cyclic finite group of order $n$ and $q$ is the least prime
divisor of $n$. We proved that 
$\psi(G)\leq\frac 7{11}\psi(C_n)$ and $\psi(G)<\frac 1{q-1}\psi(C_n)$.
The first result is best
possible, since for each $n=4k$, $k$ odd, there exists a group $G$ of order $n$
satisfying $\psi(G)=\frac 7{11}\psi(C_n)$ and the second result 
implies that if $G$ is of odd order, then  $\psi(G)<\frac 12\psi(C_n)$. 
Our results
improve the inequality  $\psi(G)<\psi(C_n)$ obtained by H. Amiri, 
S.M. Jafarian Amiri and I.M. Isaacs in 2009, as well as other results
obtained by S.M. Jafarian Amiri and M. Amiri in 2014 and by R. Shen, G. Chen 
and C. Wu in 2015.
Furthermore, we obtained some $\psi(G)$-based sufficient conditions for the 
solvability of $G$.

\end{abstract}


\end{frontmatter}

\linenumbers

\section{Introduction}


The problem of detecting structural properties of a periodic group by looking at
element orders has been considered by various authors, from many different points of view.
For example, if we denote by $\omega(G)$ the set of the orders of all the elements of $G$, there are many new and old results as well as many open questions concerning $\omega(G)$ (see for example [9]). In [1] H. Amiri, S.M. Jafarian Amiri and I.M. Isaacs introduced the function $\psi(G)$,  which denotes the sum of  element  orders  of  a  finite group $G$, and proved that if $G$ is a non-cyclic group of order $n$ then $\psi(G) < \psi(C_n)$, where $C_n$ denote the cyclic group of order $n$. Recently S.M. Jafarian Amiri and M. Amiri in [5] (see also [2] and [4]) and R. Shen, G. Chen and C. Wu in [11] studied finite groups $G$ of order $n$ with the second largest value of $\psi(G)$, and obtained information about the structure of $G$ if $n = p_1^{\alpha_1}\cdots p_t^{\alpha_t}$ , $p_1 < \cdots < p_t$, in the case $\alpha_1 > 1$.
 Products of element orders of a finite group $G$ and some other functions on the orders of the elements of $G$ have been recently studied by M. Garonzi and M. Patassini in [3].
 
 In this paper we continue the study of the function $\psi(G)$.

Our main result is the following theorem: 
\bigskip

{\bf Theorem 1}
\textit{If $G$ is a non-cyclic finite group of order $n$, then
$$\psi(G)\leq\frac 7{11}\psi(C_n).$$}

\bigskip

This upper bound is best possible, since as shown in the following
proposition, for each
$n=4k$, $k$ odd, there exists a group of order $n$ satisfying 
$\psi(G)=\frac 7{11}\psi(C_n)$.

\bigskip 

{\bf Proposition 2} \textit {Let $k$ be an odd integer  and let $n=4k$.
Then
$$\psi(C_n)= 11\psi(C_k)\qquad \psi(C_{2k}\times C_2)=7\psi(C_k)$$
and hence
$$\psi(C_{2k}\times C_2)=\frac 7{11}\psi(C_n).$$}
\bigskip

In particular, in view of Theorem 1, it follows by Proposition 2 that if
$n=4k$ with an odd $k$, then the group $G=C_{2k}\times C_2$ has the maximal sum
of element orders among non-cyclic groups of order $n$ (see also [11], Theorem 1.1) .    

 We also proved the following result, which improves Theorem 1 for groups of odd order.
\bigskip

{\bf Theorem 3} \textit{Let $G$ be a non-cyclic finite group of order $n$ and let $q$
be the smallest prime divisor of $n$. Then:
$$\psi(G)<\frac 1{q-1}\psi(C_n).$$}

\bigskip

Indeed, this theorem implies the following corollary:
\bigskip

{\bf Corollary 4}  \textit {Let $G$ be a non-cyclic finite group of odd order $n$.
Then
$$\psi(G)<\frac 12\psi(C_n).$$}

\bigskip

An important ingredient in our proofs is Corollary B of [1], which states
that if $P$ is a cyclic normal Sylow $p$-subgroup of a finite group $G$, then
$$\psi(G)\leq \psi(P)\psi(G/P),$$
with equality if and only if $P$ is central in $G$ (see Proposition 2.10). 
Another important ingredient is our
Lemma 2.1, where we proved that if $p$ and $q$ are the largest and the smallest
divisors  of an integer $n$, respectively, then the Euler's function $\varphi(n)$
satisfies the following inequality:
$$\varphi(n)\geq \frac {q-1}pn.$$

We also mention the following almost trivial upper bound for the value of $\psi(G)$ 
for non-cyclic groups $G$.

\bigskip
 {\bf Proposition 5} \textit{Let $G$ be a non-cyclic finite group of order $n$ and let $q$
be the smallest prime divisor of $n$. Then:
$$\psi(G)\leq \frac {(n-1)n}q+1<\frac {n^2}q.$$}
\medskip
Proof. 
Since $G$ is non-cyclic, it follows that $o(x)\leq n/q$ for each $x\in
G$.
But $o(1)=1$, so $\psi(G)\leq (n-1)(n/q)+1<n^2/q$, as required.
\qed

\bigskip

Notice, however that $\psi(S_3)=13>\frac 12(q-1)\psi(C_6)=\frac {21}2$. 
This observation raises the following question: what can we say about groups of order
$n$ satisfying
$\psi(G)\geq (1/2(q-1))\psi(C_n)$?
A partial answer can be found in our next theorem.

\bigskip
{\bf Theorem 6} \textit{Let $G$ be  a finite group of order $n$  and let $q$ and $p$
be the smallest and the largest prime divisors of $n$, respectively. Suppose that
$G$ satisfies
$$\psi(G)\geq \frac{1}{2(q-1)}\psi(C_n).$$
Then $G$ is solvable, the Sylow $p$-subgroups of $G$ contain a cyclic subgroup of index
$p$ and one of the following statements holds:}
\begin{enumerate}

\item \textit{ The Sylow $p$-subgroup $P$ of $G$ is cyclic and normal in $G$};

\item \textit{The Sylow $q$-subgroups of $G$ are cyclic, $G$ is $q$-nilpotent and
$G^{\prime\prime}\leq Z(G)$};

\item \textit{The Sylow $p$-subgroups of $G$ are cyclic, $G$ is $p$-nilpotent and
$G^{\prime\prime}\leq Z(G)$.}
\end{enumerate}

Theorem 6 implies the following corollaries. We use $p$ and $q$ as defined in
Theorem 6.

\bigskip

{\bf Corollary 7} \textit{The conclusions of Theorem 6 hold if $G$ satisfies
$$\psi(G)\geq \frac 1q\psi(C_n).$$}
\medskip
Proof. 
Since $q\geq 2$, it follows that $q\leq 2(q-1)$.
\qed

\bigskip

{\bf Corollary 8} \textit{The conclusions of Theorem 6 hold if $G$ is a group of 
odd order satisfying
$$\psi(G)\geq \frac 1{q+1}\psi(C_n).$$}
\medskip
Proof.  Since $q\geq 3$, it follows that $q+1\leq 2(q-1)$.
\qed

\bigskip
  
{\bf Corollary 9} \textit{If either $G$ is non-solvable or a Sylow $p$-subgroup of
$G$ contains no cyclic subgroup of index $p$, then
$$\psi(G)<\frac 1{2(q-1)}\psi(C_n)\leq \frac 1q\psi(C_n).$$}

\bigskip

Our next result is another $\psi(G)$-based sufficient condition for the solvability
of $G$.  

\bigskip

{\bf Theorem 10} \textit{Let $G$ be a finite group of order $n$ satisfying
$$\psi(G)\geq \frac 35n\varphi(n).$$
Then $G$ is solvable and $G^{\prime\prime}\leq Z(G)$.}

\bigskip

This condition is certainly not necessary for the solvability of $G$.
For example, for $n=8$ we have 
$$\psi(C_2\times C_2\times C_2)=15<\frac 35\cdot 8\cdot 4=\frac 35n\varphi(n).$$
On the other hand, for $n=60$, the simple group $A_5$ satisfies 
$\psi(A_5)=211>\frac 15n\phi(n)=192$. 
 
In the proof of Theorem 10 we apply the following result of Ramanujan
(see [8], page 46):
if $q_1= 2, q_2,
\cdots , q_n, \cdots$ is the increasing sequence of all primes, then            
$$\prod_{i =1,\cdots,\infty} \frac{q_i^2+1}{q_i^2-1} = \frac 52.$$

Our final result deals with groups of order $n$ which satisfy $\psi(G)\geq \frac 1q
n\varphi(n)$.

\bigskip
{ \bf Theorem 11} \textit{ Let $G$ be  a finite group of order $n$  and let $q$ and $p$
be the smallest and the largest prime divisors of $n$, respectively. Suppose that
$G$ satisfies
$$\psi(G)\geq \frac 1q n\varphi(n).$$
Then either $G$ has a normal cyclic Sylow $p$-subgroup or it 
is a solvable group with
a cyclic maximal subgroup of index either $p$ or $p+1$.}

\section*{ Acknowledgment}

We are grateful to our colleague Guy Moshkowich for constructing the table 
of the values of the function $\psi(G)$ for all groups $G$ of order less than $128$. 
This table,
which he constructed using the GAP System of Computational Group Theory, provided useful
information for our research.

This work was supported by the "National Group for Algebraic and Geometric Structures, and their Applications" (GNSAGA - INDAM), Italy.

The first author is grateful to the Department of Mathematics of the University of Salerno for its hospitality and support,
while this investigation was carried out.

\section{Preliminary results}


First we determine a lower bound for  $\varphi(n)$.
\bigskip

 {\bf Lemma 2.1} \textit{Let $n$ be a positive integer larger than $1$,
with the largest prime divisor $p$ and the smallest prime divisor $q$.
Then
$$\varphi(n)\geq \frac {q-1}pn.$$}

\medskip
Proof. Let $n=p_1^{r_1}p_2^{r_2}\dotsb p_k^{r_k}$, where the $p_i$'s are primes,
the $r_i$'s are positive integers and $p=p_1>p_2>\dots >p_k=q$. Our proof is by
induction on
$k$.

If $k=1$, then $n=p^{r_1}$ and
$$\varphi(n)=\varphi(p^{r_1})=\frac {p-1}{p}p^{r_1}=\frac {p-1}{p}n,$$ 
as required.

Suppose now that $k>1$ and that the lemma holds for all integers with less 
than $k$ distinct prime
divisors. Set $m=p_2^{r_2}\dotsb p_k^{r_k}$. Then by induction $\varphi(m)\geq
\frac {p_k-1}{p_2}m$ and
$$\varphi(n)=\varphi(p_1^{r_1})\varphi(m)\geq \frac{p_1-1}{p_1}p_1^{r_1}\frac
{p_k-1}{p_2}m
\geq \frac{p_1-1}{p_1}\frac{p_k-1}{p_1-1}n=\frac{p_k-1}{p_1}n=\frac{q-1}pn,$$
as required. The proof is now complete.
\qed

\bigskip

Our next aim is to find a convenient formula for $\psi(G)$ when $G=P\rtimes F$,
$P$ is a cyclic $p$-group for some prime $p$, $|F|>1$ and $(p,|F|)=1$.

\bigskip

{\bf Lemma 2.2} \textit{Let $G$ be a finite group satisfying $G=P\rtimes F$, where $P$
is a cyclic $p$-group for some prime $p$, $|F|>1$ and $(p,|F|)=1$. Then the following
statements hold.}
\begin{enumerate}

\item [$(1)$] \textit {Each element of $F$ acts on $P$ either trivially or fixed-point-freely.}
\item [$(2)$]  \textit{If $x\in F$, $o(x)=m$ and $u\in P$, then $m$ is the least 
positive integer satisfying $(ux)^m\in P$.}
\item [$(3)$]  \textit{If $u\in P$ and $x\in C_F(P)$,  then $o(ux)=o(u)o(x)$.}
\item [$(4)$] \textit{If $u\in P$ and $x\in F\setminus C_F(P)$, then $o(ux)=o(x)$.}
\item [$(5)$]  \textit {Let $Z=C_F(P)$. Then
$$\psi(G)=\psi(P)\psi(Z) + |P|\psi(F\setminus Z)<\psi(P)\psi(Z) + |P|\psi(F).$$}
\end{enumerate}
\medskip
Proof.

(1)  Suppose that $x\in F$ acts trivially on $u\in P\setminus \{1\}$. Then $x$ acts
trivially on $\Omega_1(P)$ and hence it  acts trivially on $P$ (see [6], Theorem
5.4.2).
The claim follows from this remark.

(2) Since $P\triangleleft G$, it follows that if $n$ is a positive integer,
then $(ux)^n=v_nx^n $ for some $v_n\in P$. As
$P\cap F=\{1\}$, it follows that $(ux)^n\in P$ if and only if $m$ divides $n$. The
claim
follows.

(3) Trivially holds.

(4) Suppose that $o(x)=m$. By (2) $(ux)^m\in P$ and hence
$1=[(ux)^m,ux]=[(ux)^m,x]$. Since $x\in F\setminus C_F(P)$, it follows by (1) that
$(ux)^m=1$ and by (2) $o(ux)=m=o(x)$, as required.

(5) It follows by (3)  that $\psi(PZ)=\psi(P)\psi(Z)$
and by (4) that $\psi(G\setminus (PZ))=|P|\psi(F\setminus Z)$. Therefore
$$\psi(G)=\psi(P)\psi(Z) + |P|\psi(F\setminus Z)<\psi(P)\psi(Z) + |P|\psi(F).$$

\qed

\bigskip

We also need information concerning finite groups with a cyclic maximal
subgroup. First we mention the following related result of Herstein [7].
\bigskip

{\bf Proposition 2.3} (Herstein) \textit{If $G$ is a
finite group with an abelian maximal
subgroup, then $G$ is solvable.}

\bigskip
Using this result we proved the following proposition, which is of independent
interest.

\bigskip
{\bf Proposition 2.4} \textit{Let $G$ be a finite group with a cyclic maximal
subgroup $C$. Then $G$ is solvable and $G^{\prime\prime}\leq Z(G)$.}

\medskip
Proof. The group $G$ is solvable by Proposition 2.3.
If $G^{\prime}\leq C$, then $G^{\prime\prime}=1\leq Z(G)$ as required.
Otherwise $G= G^{\prime}C$ and $G^{\prime\prime}\leq C$, since otherwise
$G=G^{\prime\prime}C$ and $G^{\prime}\leq G^{\prime\prime}$, a contradiction
since $G$ is solvable and $G^{\prime}\neq 1$. Hence $G^{\prime\prime}$ is cyclic and
$G/C_G(G^{\prime\prime})$ is abelian. Consequently $G^{\prime}$ and $C$
are both subgroups of $C_G(G^{\prime\prime})$, yielding
again $G^{\prime\prime}\leq Z(G)$, as required. 
\qed

\bigskip

Another important and useful result is the following proposition.

\bigskip

{\bf Proposition 2.5} \textit{Let $G$ be a finite group and suppose that there exists
$x\in G$ such that
$$[G:\langle x\rangle]<2p,$$
where $p$ is the maximal prime divisor of $|G|$. Then one of the following holds:}

\textit{(i) $G$ has a normal cyclic Sylow $p$-subgroup,}

\textit{(ii) $G$ is solvable and $\langle x\rangle$ is a maximal subgroup of $G$ 
of index either $p$ or $p+1$.}

\medskip

Proof. First suppose that $p$ divides $[G:\langle x\rangle]$. Since
$[G:\langle x\rangle]$ divides $|G|$, our assumption implies that
$[G:\langle x\rangle]=p$ and $G$ is solvable by Proposition 2.4. Thus $G$
satisfies (ii).

Now assume that $p$ does not divide  $[G:\langle x\rangle]$. Then $\langle x\rangle$
contains a cyclic Sylow $p$-subgroup $P$ of $G$. If $P$ is normal in $G$, then (i)
holds. So suppose, finally, that $P$ is not normal in $G$. Since
$\langle x\rangle \leq N_G(P)$, it follows from our assumptions that
$[G:N_G(P)]<2p$. Since $P$ is not normal in $G$, this implies 
that $[G:N_G(P)]=p+1$ and that $N_G(P)$
is a maximal subgroup of $G$. 
But
$$[N_G(P):\langle x\rangle]=\frac {[G:\langle x\rangle]}{[G:N_G(P)]}
<\frac {2p}{p+1}<2,$$
so $N_G(P)=\langle x\rangle$ and $\langle x\rangle$
is a cyclic maximal subgroup of $G$ of index $p+1$. By Proposition 2.4 
$G$ is solvable, and hence it satisfies (ii). 
The proof of the proposition is complete.
\qed

\bigskip

We also need the following related result.

\bigskip

{\bf Proposition 2.6} \textit{The following statements hold.}
\begin{enumerate}
\item [$(1)$] \textit{If $G$ is a finite $2$-group with a cyclic subgroup of index $4$,
then $G^{\prime\prime}\leq Z(G)$.}
\item [$(2)$] \textit{ If $G$ is a finite group of order $2^{\alpha}3^{\beta}$ with a cyclic
subgroup of index less than $6$, then $G^{\prime\prime}\leq Z(G)$.}
\end{enumerate}
\medskip

Proof.
(1) Let $\langle a \rangle$ be of index $4$ in $G$ and let $M$ be a maximal subgroup
of $G$ containing $\langle a \rangle$. If $M=\langle a \rangle$, then
the result follows by Proposition 2.4. So assume that $M>\langle a \rangle$.
Then $[G:M]=2$, implying that $M$ is normal in $G$,
$G^{\prime} \leq M$ and $G^{\prime\prime} \leq M^{\prime}$.
Moreover, $M$ has a maximal
cyclic subgroup and therefore by Theorem 5.3.4 in [10], either $M$ is abelian,
or
$M^{\prime}$
has order 2, or $M$ is dihedral, semidihedral or generalized quaternion.

If $M$ is abelian then $G^{\prime\prime}=1$ and if $|M^{\prime}| \leq 2$,
then $G^{\prime\prime} \leq M^{\prime}$ implies that
$G^{\prime\prime}$ is a normal subgroup of
$G$ of order at most $2$, thus $G^{\prime\prime} \leq Z(G)$, as required.

Now suppose that $M$ is either dihedral, or semidihedral or
generalized quaternion.
Then there exists $x\in M$ such that
$a^x = a^{-1}a^{\gamma2^{n-1}}$, where $o(a) = 2^n$, $\gamma \in \{0, 1\}$,
$o(x)\in \{2, 4 \}$ and $x^2\in Z(M)$.

Write $G = M \langle y \rangle$.

If $a^y \in \langle a \rangle$, then $\langle a \rangle$ is normal in $G$,
thus
$G^{\prime} \leq \langle a \rangle$ since $|G/ \langle a \rangle| = 4$.
Hence $G^{\prime}$ is abelian and $G^{\prime\prime}=1$, as required.

Suppose, finally, that  $a^y \notin \langle a \rangle$.
Then $a^y = a^\delta x$, where $\delta$ is an integer.
We have
$$(a^2)^y = a^\delta xa^\delta x = a^\delta x^2 a^{-\delta}a^{\gamma 2^{n-1}\delta
}  = x^2 a^{\gamma 2^{n-1}\delta } $$
and $(a^y)^4 =
(a^{\gamma 2^{n-1} \delta} x^2 )^2 = x^4 = 1$. Hence $o(a) = o(a^y) = 4$,
$|M|=8$ and
$M^{\prime} \leq \langle a^2 \rangle$. Thus $G^{\prime\prime} \leq
M^{\prime}$
has 
order at most $2$
and hence it is contained in $Z(G)$, as required.

(2) Let $\langle a \rangle$ be a cyclic subgroup of $G$ of index less than
$ 6$.

If $|G : \langle a \rangle| = 2$ or $|G : \langle a \rangle| = 3$, then the
result
follows from Proposition 2.4.

Suppose, finally, that $|G : \langle a \rangle| = 4$.
If $\beta =0$ then the result follows by (1), and if
$\langle a \rangle$ is maximal in $G$,
then the result follows again from Proposition 2.4.

So suppose that $\beta >0$ and that $\langle a \rangle$ is not maximal in
$G$.
Then by  Proposition 2.5 $G$ has a normal
cyclic Sylow $3$-subgroup $P$. Thus $G = P \rtimes D$ ,
where
$|D| = 2^\alpha$. Obviously $P \leq \langle a \rangle$ and $D \simeq G/P$ has
the
cyclic subgroup $\langle a \rangle/P$ of index $4$. Hence by (1),
$D^{\prime\prime}\leq Z(D)$. Now we have $G = PD$, $G^{\prime} \leq C_G(P)$
and $G^{\prime} =D^{\prime}[P,D]$, which implies that
$$G^{\prime\prime} = D^{\prime\prime} \leq Z(D)\cap C_G(P) \leq Z(G),$$ as
required.
\qed

\bigskip

We also state the result of Ramanujan (see [8], page 46), which was
mentioned in the introduction.

\bigskip

{\bf Proposition 2.7 }(Ramanujan)
\textit{If $q_1= 2, q_2,
\cdots , q_n, \cdots$ is the increasing sequence of all primes, then
$$\prod_{i =1,\cdots,\infty} \frac{q_i^2+1}{q_i^2-1} = \frac 52.$$}

\bigskip

This proposition implies the following lemma.

\bigskip

{\bf Lemma 2.8} \textit{Let $p_2,p_3,\dots,p_s$ be primes satisfying
$p_2<p_3<\dotsb <p_s$. If $p_2>3$ then
$$\prod_{i = 2}^{s} \frac{p_i^2-1}{p_i^2+1} >\frac 56.$$}

Proof. If $p_2>3$, then 
Proposition 2.7 implies that
$$\frac{2^2+1}{2^2-1}\frac{3^2+1}{3^2-1}\prod_{i = 2,\cdots,s}
\frac{p_i^2+1}{p_i^2-1} < \frac 52.$$
Thus $\prod_{i = 2,\cdots, s} \frac{p_i^2+1}{p_i^2-1}<\frac 65$,
yielding
$$\prod_{i= 2,\cdots,s} \frac{p_i^2-1}{p_i^2+1} >\frac 56,$$
as required.
\qed

\bigskip

We shall also need some basic facts about $\psi(C_n)$.

\bigskip

{\bf Lemma 2.9}
\begin{enumerate}
\item [$(1)$] \textit{If $P$ is a cyclic group of order $p^r$ for some prime $p$, then
$$\psi(P)=\frac {p^{2r+1}+1}{p+1}=\frac {p|P|^2+1}{p+1}.$$}

\item [$(2)$] \textit{Let $p_1<p_2<\dots <p_t=p$ be the prime divisors of $n$ and denote
the corresponding Sylow subgroups  of $C_n$ by $P_1,P_2,\dots,P_t$. 
Then 
$$\psi(C_n)=\prod_{i=1}^t \psi(P_i)\geq \frac 2{p+1}n^2.$$}
\end{enumerate}

\bigskip
Proof. 

(1) $\psi(P)=1+p\varphi(p)+p^2\varphi(p^2)+\dotsb +p^r\varphi(p^r)=
\frac {p^{2r+1}+1}{p+1}=\frac {p|P|^2+1}{p+1}$.

(2) Since $C_n=P_1\times P_2\times \dotsb\times P_t$, it follows by Lemma 2.2(3)
that $\psi(C_n)=\prod_{i=1}^t \psi(P_i)$. 
Since $p_{i+1}\geq p_i+1$ for all $i$ and $p_1\geq 2$, it follows by (1) that
$$\psi(C_n)=\prod_{i=1}^t \frac {p_i|P_i|^2+1}{p_i+1}
>\prod_{i=1}^t \frac {p_i}{p_i+1}|P_i|^2\geq \frac 2{p+1}n^2.$$
\qed

\bigskip

In most cases, we shall apply the results of Lemma 2.9 without  reference.
 
Finally, we state  Corollary B from [1].

{\bf Proposition 2.10} \textit{If $P$ is a cyclic normal Sylow $p$-subgroup 
of a finite group $G$, then
$$\psi(G)\leq \psi(P)\psi(G/P),$$
with equality if and only if $P$ is central in $G$.}

\bigskip

\section{Proofs of the main results.}


\bigskip

Since we are using the result of Theorem 3 for the proof of Theorem 1, we
shall prove Theorem 3 first. The proof of Proposition 2 will follow that of
Theorem 1.
\bigskip

\textit{Proof of Theorem 3. }We need to prove that if $\psi(G)\geq\frac
1{q-1}\psi(C_n)$,
then $G\cong C_n$.

Clearly $\psi(C_n)>n\varphi(n)$ and by Lemma 2.1 $\varphi(n)\geq (q-1)n/p$, where
$p$ denotes the largest prime divisor of $n$. Hence by our assumptions
$\psi(G)>\frac {n(q-1)n}{(q-1)p}=n^2/p$, which implies that there exists $x\in G$
with $o(x)>n/p$. Thus $[G:\langle x\rangle]<p$ and $\langle x\rangle$ contains
a Sylow $p$-subgroup $P$ of $G$. Since $\langle x\rangle\leq N_G(P)$, it follows that
$P$ is a cyclic normal subgroup of $G$ and Proposition 2.10 implies that
$$\psi(P)\psi(G/P)\geq \psi(G)\geq \frac 1{q-1}\psi(C_{p^r})\psi(C_{n/p^r}),$$
where $p^r=|P|$. Since $P\cong C_{p^r}$, cancellation yields
$$\psi(G/P)\geq \frac 1{q-1}\psi(C_{n/p^r}).$$

If $n=p^r$, $p$ a prime, then the existence of $x\in G$ satisfying $o(x)>n/p$
implies that $o(x)=n$ and $G$ is cyclic, as required. So we may assume that
$n$ is divisible by exactly $k$ different primes with $k>1$. Applying induction
with respect to $k$, we may assume that the theorem holds for groups of order
which has less than
$k$ distinct prime divisors. Since $|G/P|$ has $k-1$ distinct prime divisors and $G/P$
satisfies our assumptions, it follows that $G/P$ is cyclic and $G=P\rtimes F$,
with $F\cong G/P$ and $F\neq 1$. Notice that $n=|P||F|$, $P$ and $F$ are
both cyclic  and $(|P|,|F|)=1$. Hence
$\psi(C_n)=\psi(P)\psi(F)$.

If $C_F(P)=F$, then $G=P\times F$ and $G$ is cyclic, as required.

So it suffices to prove that if $C_F(P)=Z<F$, then  $\psi(G)<(1/(q-1))\psi(C_n)$,
contrary to our assumptions. It follows by Lemma 2.2(5) that
$$\psi(G)=\psi(P)\psi(Z) + |P|\psi(F\setminus Z)<\psi(P)\psi(Z) + |P|\psi(F).$$
Hence
$$\psi(G)<\psi(P)\psi(F)(\frac {\psi(Z)}{\psi(F)}+\frac {|P|}{\psi(P)})
= \psi(C_n)(\frac {\psi(Z)}{\psi(F)}+\frac {|P|}{\psi(P)}).$$
Notice first that since $P$ is a cyclic $p$-group,
we have
$$\frac {|P|}{\psi(P)}=\frac {|P|(p+1)}{p|P|^2+1}<\frac {p+1}{p|P|}\leq
\frac {p+1}{p^2}<\frac {p+1}{p^2-1}= \frac 1{p-1}\leq \frac 1q.$$
Next notice that $Z$ is a proper subgroup of the cyclic group $F$ and $\psi(F)$
is a product of $\psi(S)$, with $S$ running over the Sylow subgroups of $F$.
Since also $\psi(Z)$ is a similar product, and since at least one Sylow
subgroup of $Z$, say Sylow $r$-subgroup $R_Z$, is properly contained in the Sylow
$r$-subgroup $R_F$ of $F$ of order $r^s$, it follows that
$$\frac {\psi(Z)}{\psi(F)}\leq \frac {\psi(R_Z)} {\psi(R_F)}\leq \frac
{r^{2(s-1)+1}+1}
{r^{2s+1}+1}.$$
Since $r\geq  q$ and $s\geq 1$, we get
$$\frac {\psi(Z)}{\psi(F)}\leq \frac {r^{2s-1}+1}{r^{2s+1}+1}
<\frac 1{q(q-1)}.$$
Therefore
$$\psi(G)<\psi(C_n)(\frac {\psi(Z)}{\psi(F)}+\frac {|P|}{\psi(P)})<
\psi(C_n)(\frac 1{q(q-1)}+\frac 1q)=\psi(C_n)\frac 1{q-1},$$
a contradiction.

The proof is now complete.
\qed

\bigskip

We continue with the proof of Theorem 1.

\bigskip

\textit{Proof of Theorem 1.}
 Throughout this proof $G$ denotes a non-cyclic finite group of order $n$
satisfying
$$\psi(G)>\frac 7{11}\psi(C_n).$$ 
Let $p_1<p_2<\dots <p_t=p$ be the prime divisors of $n$ and denote
the corresponding Sylow subgroups  of $C_n$ by $P_1,P_2,\dots,P_t$.
By Lemma 2.9 
$\psi(C_n)\geq\frac 2{p+1}n^2$,
so our assumptions imply that $G$ satisfies
$$\psi(G)>\frac 7{11}\psi(C_n)\geq  \frac {14}{11(p+1)}n^2.$$
Our aim is to reach a contradiction. Our proof is by induction on the size of $p$.
 
By our assumptions there exists $x\in G$
with
$o(x)>\frac {14}{11(p+1)}n$,
which implies that  
$$[G:\langle x\rangle]<\frac {11(p+1)}{14}.$$ 

Suppose, first, that $p=2$. Then $G$ is a $2$-group and 
$[G:\langle x\rangle]<\frac {33}{14}$. 
Thus $[G:\langle x\rangle]=2$, $n\geq 4$
and $n^2\geq 16$, implying that 
$$\psi(G)\leq \psi(C_{n/2})+(\frac n2)^2=\frac {2(n/2)^2+1}3+\frac {n^2}4
=\frac 5{12}n^2+\frac 13 \leq
(\frac 7{11})(\frac {2n^2+1}3)=(\frac 7{11})\psi(C_n),$$
a contradiction.

Next assume that $p=3$ and $[G:\langle x\rangle]<\frac {44}{14}$. 
If $G$ is a $3$-group, then Theorem 3 yields
$\psi(G)<\frac 12\psi(C_n)<\frac 7{11}\psi(C_n)$, a contradiction.
So we may assume that $n=2^a3^b$ for some positive integers $a$ and $b$.
Since $[G:\langle x\rangle]<\frac {44}{14}$, it follows that
$[G:\langle x\rangle]\leq 3$. Hence either
$[G:\langle x\rangle]=2$ or $[G:\langle x\rangle]=3$.
Notice for later reference that
if $n=2^a3^b$, then
$$
\frac 7{11}\psi(C_n)$$ $$=\frac 7{11}\psi(C_{2^a})\psi(C_{3^b})
=\frac 7{11}(\frac {2^{2a+1}+1}3)(\frac {3^{2b+1}+1}4)$$
$$=\frac 7{22}2^{2a}3^{2b}+\frac 7{66}2^{2a}+\frac 7{44}3^{2b}+\frac 7{132}.
$$

Suppose first that $[G:\langle x\rangle]=2$.
Then $\langle x\rangle$ contains a cyclic Sylow $3$-subgroup $P$ of $G$ and
since $\langle x\rangle\leq C_G(P)$, it follows that $P$ is normal in $G$.

If there exists $y\in G\setminus \langle x\rangle $ with $[G:\langle y\rangle]=2$,
then $y\in C_G(P)$ and hence $P\leq Z(G)$. Thus $G=P\times Q$, where $Q$ is a
non-cyclic Sylow $2$-subgroup of $G$ and it follows by our result for $p=2$
that
$$\psi(G)=\psi(P)\psi(Q)\leq \psi(P)(\frac 7{11})\psi(C_{|Q|})
=(\frac 7{11})\psi(C_n),$$
a contradiction.

So suppose that   $o(y)\leq \frac n3$ for all $y\in G\setminus \langle x\rangle$. 
Since $a\geq 1$ and $b\geq 1$,
we obtain the following final contradiction with respect to $[G:\langle x\rangle]=2$:
$$ 
\psi(G) \leq \psi(C_{\frac n2})+(\frac n2)(\frac n3)
=\psi(C_{2^{a-1}})\psi(C_{3^b})+\frac {n^2}6
=(\frac {2^{2a-1}+1}3)(\frac{3^{2b+1}+1}4)+\frac {n^2}6$$ $$=(\frac 16)(\frac 34)2^{2a}3^{2b}+\frac {2^{2a}}{24}
+\frac {3^{2b}}4+\frac 1{12}+\frac {2^{2a}3^{2b}}6
=\frac 7{24}2^{2a}3^{2b} +\frac {2^{2a}}{24}+\frac {3^{2b}}4+\frac 1{12}$$ $$=\frac 7{22}2^{2a}3^{2b}+\frac 7{66}2^{2a}+\frac 7{44}3^{2b}+\frac 7{132}$$
$$+(\frac 7{24}-\frac 7{22})2^{2a}3^{2b}+(\frac 1{24}-\frac 7{66})2^{2a}
+(\frac 14-\frac 7{44})3^{2b}+(\frac {11}{132}-\frac 7{132})$$
$$<\frac 7{11}\psi(C_n)-\frac 7{264}2^{2a}3^{2b}+\frac 1{11}3^{2b}
+\frac 4{132}$$
$$\leq \frac 7{11}\psi(C_n)-\frac 7{66}3^{2b}+\frac 6{66}3^{2b}
+\frac 2{66}<\frac7{11}\psi(C_n). $$

We approach now the second possibility: $p=3$ and $[G:\langle x\rangle]=3$.
By the previous arguments we may assume that
no element of $G$ is of order $\frac n2$  and hence
$o(y)\leq \frac n3$ for all $y\in G$.
By considering elements of $G$ belonging to  $\langle x\rangle$ and
those outside it, and recalling that $b\geq 1$, we obtain the following
contradiction:
$$
\psi(G)\leq \psi(C_{2^a})\psi(C_{3^{b-1}}) +2(\frac n3)^2
=(\frac {2^{2a+1}+1}3)(\frac {3^{2b-1}+1}4)+\frac 29n^2$$
$$=\frac 1{18}2^{2a}3^{2b}+\frac 162^{2a}
+\frac 1{36}3^{2b}+\frac 1{12}+\frac 292^{2a}3^{2b}
=\frac 5{18}2^{2a}3^{2b}+\frac 162^{2a}+\frac 1{36}3^{2b}+\frac 1{12}$$
$$=\frac 7{22}2^{2a}3^{2b}+\frac 7{66}2^{2a}+\frac 7{44}3^{2b}+\frac 7{132}$$
$$+(\frac 5{18}-\frac 7{22})2^{2a}3^{2b}+(\frac 16-\frac 7{66})2^{2a}
+(\frac 1{36}-\frac 7{44})3^{2b}+(\frac {11}{132}-\frac 7{132})$$
$$<\frac 7{11}\psi(C_n)-\frac 4{99}2^{2a}3^{2b}+\frac 2{33}2^{2a}+\frac 1{33}$$
$$\leq\frac 7{11}\psi(C_n)-\frac {36}{99}2^{2a}+\frac 6{99}2^{2a}+\frac 1{33}
<\frac 7{11}\psi(C_n).$$

So assume, finally, that $p>3$ and that the theorem holds for smaller values
of $p$.  Then
$$ [G:\langle x\rangle]<\frac {11(p+1)}{14}\leq p$$
and $\langle x\rangle$ contains
a cyclic Sylow $p$-subgroup $P$ of $G$. Since $\langle x\rangle\leq N_G(P)$, it follows that
$P$ is a cyclic normal subgroup of $G$ and Proposition 2.10 implies that 
$$\psi(P)\psi(G/P)\geq \psi(G)> \frac 7{11}\psi(C_{p^r})\psi(C_{n/p^r}),$$
where $p^r=|P|$. Since $P\cong C_{p^r}$, cancellation yields
$$\psi(G/P)> \frac 7{11}\psi(C_{n/p^r}).$$
Since the maximal prime dividing $n/p^r$ is smaller than $p$, our
induction hypothesis implies that $G/P$ is cyclic and
$G=P\rtimes F,$
with $F\cong G/P$ and $F\neq 1$. Notice that $n=|P||F|$, with both
$P$ and $F$ being cyclic,  and $(|P|,|F|)=1$. Hence
$\psi(C_n)=\psi(P)\psi(F)$.

If $C_F(P)=F$, then $G=P\times F$ and $G$ is a cyclic, a contradiction.
So suppose that $C_F(P)=Z<F$. Lemma 2.2(5) then implies that

$$\psi(G)<\psi(P)\psi(F)(\frac {\psi(Z)}{\psi(F)}+\frac {|P|}{\psi(P)})
= \psi(C_n)(\frac {\psi(Z)}{\psi(F)}+\frac {|P|}{\psi(P)}).$$
Notice first that since $P$ is a cyclic $p$-group and $p>3$,
we have
$$\frac {|P|}{\psi(P)}=\frac {|P|(p+1)}{p|P|^2+1}<\frac {p+1}{p|P|}\leq
\frac {p+1}{p^2}\leq \frac 6{25}<\frac 14.$$
Next notice that $Z$ is a proper subgroup of the cyclic group $F$ and $\psi(F)$
is a product of $\psi(S)$, with $S$ running over the Sylow subgroups of $F$.
Since also $\psi(Z)$ is a similar product, and since at least one Sylow
subgroup of $Z$, say Sylow $r$-subgroup $R_Z$, is properly contained in the Sylow
$r$-subgroup $R_F$ of $F$ of order $r^s$, it follows that
$$\frac {\psi(Z)}{\psi(F)}\leq \frac {r^{2(s-1)+1}+1}{r^{2s+1}+1}.$$
But $r\geq  2$ and $s\geq 1$, so
$$\frac {r^{2(s-1)+1}+1}{r^{2s+1}+1}\leq \frac 1{r+1}, $$
since this inequality is equivalent to $1\leq r^{2s-2}(r^2-r-1)$,
which is true. Hence
$$\frac {\psi(Z)}{\psi(F)}\leq \frac 1{r+1}\leq \frac 13,$$
and
$$\psi(G)<\psi(C_n)(\frac {\psi(Z)}{\psi(F)}+\frac {|P|}{\psi(P)})<
\psi(C_n)(\frac 13+\frac 14)=\psi(C_n)\frac 7{12}<\psi(C_n)\frac 7{11},$$
a final contradiction.

The proof is now complete.
\qed

\bigskip

Next we prove Proposition 2.
\bigskip

\textit{Proof of Proposition 2}. We start with the proof of the first equality:
$$\psi(C_n)= \psi(C_{4k})=\psi(C_4)\psi(C_k)=\frac
{32+1}{2+1}\psi(C_k)=11\psi(C_k).$$
Next we prove the second equality:
$$\psi(C_{2k}\times C_2) =\psi(C_{k} \times C_{2} \times C_{2}) =  \psi(C_{k})\psi(C_{2} \times C_2)= 7\psi(C_{k}).
$$
The claim follows.
\qed

\bigskip

We continue with a proof of Theorem 6.

\bigskip

\textit{Proof of Theorem 6.}
Since  $\psi(C_n)>n\varphi(n)$ and by Lemma 2.1 $\varphi(n)\geq \frac {(q-1)n}p$,
it follows  by our assumptions that
$\psi(G)>\frac {n^2}{2p}$. Hence there exists $x\in G$ such that
$o(x)>\frac n{2p}$ and 
$$[G:\langle x\rangle]< 2p.$$

First we shall prove that $G$ is solvable. Let $k$ be the number of prime
divisors of $n$. Our proof is by induction on $k$.

If $k=1$, then $G$ is a $p$-group, hence solvable, as claimed.
So assume that $k>1$ and that the claim holds for $k-1$.

If $p\mid [G:\langle x\rangle]$, then $[G:\langle x\rangle]=p$ and
$\langle x\rangle$ is a cyclic maximal subgroup of $G$. Hence by Proposition 2.3
$G$ is solvable, as claimed.

So suppose that $p\nmid [G:\langle x\rangle]$. Then $\langle x\rangle$ 
contains a cyclic Sylow $p$-subgroup $P$ of $G$. 

If $P$ is normal in $G$, then Proposition 2.10 and our assumptions imply that
$$\psi(P)\psi(G/P)\geq \psi(G)\geq \frac 1{2(q-1)}\psi(C_{|P|})\psi(C_{|G/P|}),$$
and since $\psi(P)=\psi(C_{|P|})$, it follows that
$$\psi(G/P)\geq \frac 1{2(q-1)}\psi(C_{|G/P|}).$$
Hence, by induction, $G/P$ is solvable  and so $G$ is solvable, as claimed.

Suppose, finally, that $P$ is not normal in $G$. Since $\langle x\rangle$  is a 
subgroup of $N_G(P)$, it follows that $[G:N_G(P)]<2p$ and hence 
$$[G:N_G(P)]=p+1.$$
Since
$$[N_G(P):\langle x\rangle]=\frac {[G:\langle x\rangle]}{[G:N_G(P)]}\leq \frac
{2p}{p+1}<2,$$
it follows that $N_G(P)=\langle x\rangle$ is a cyclic maximal subgroup of $G$.
Hence $G$ is solvable by Proposition 2.3 and the proof of our claim is complete.

We proceed with the proof of the theorem. As shown above, there exists $x\in G$
such that $[G:\langle x\rangle]<2p$. By Proposition 2.5 this implies that 
either $G$ has a normal cyclic Sylow $p$-subgroup or $\langle x\rangle$ is a maximal 
subgroup of $G$ of index either $p$ or $p+1$. In either case, the Sylow $p$-subgroups
of $G$ contain a cyclic subgroup of index $p$, as required.

If $G$ has a normal cyclic Sylow $p$-subgroup, then (1) holds. If 
$\langle x\rangle$ is a maximal
subgroup of $G$ of index $p$, then it contains a Sylow $q$-subgroup $Q$ of $G$.
Since $Q$ is cyclic and $q$ is the smallest prime divisor of $n$, it follows 
by Theorem 10.1.9 in [5] that $G$
is $q$-nilpotent. Moreover, by Proposition 2.4 $G^{\prime\prime}\leq Z(G)$  and 
(2) holds.

Finally, if $\langle x\rangle$ is a maximal
subgroup of $G$ of index $p+1$, then it contains a cyclic Sylow $p$-subgroup $P$ of $G$.
If $P$ is normal, then (1) holds. So suppose that $P$ is not normal in $G$.
Since $\langle x\rangle\leq N_G(P)$, it follows that $[G:N_G(P)]<2p$, which implies
that $[G:N_G(P)]=p+1$. As shown above $N_G(P)=\langle x\rangle$ and hence 
$N_G(P)=C_G(P)$. Thus $G$ is $p$-nilpotent by Burnside's theorem, and since
$\langle x\rangle$ is a cyclic maximal  subgroup of $G$, Proposition 2.4 implies
that $G^{\prime\prime}\leq Z(G)$ and (3) holds. The proof of the theorem is
complete.
\qed

\bigskip

We continue with a proof of a sufficient condition for the solvability 
of a finite group $G$.

\bigskip

\textit{Proof of Theorem 10. } Suppose that 
$$\psi(G)\geq \frac 35n\varphi(n)$$
and let $p_1$ be the maximal prime dividing $n$.
By Lemma 2.1
$\varphi(n) \geq n/p_1$, so by our assumptions $\psi(G) \geq\frac 35 n^2/p_1$.
Hence there exists
an element $x$ of $G$ with
$o(x)>\frac 35 n/p_1$ and
$$|G: \langle x \rangle| < \frac 53 p_1 < 2p_1.$$
It follows  by Proposition 2.5 that either $G$ is solvable,  or $G$
has a normal cyclic Sylow $p_1$-subgroup $P_1$.

We prove first that $G$ is solvable. Clearly we may assume that $n$ is
divisible by at least three different primes
and hence $p_1\geq 5$. We may also assume that $G$ has
a normal cyclic Sylow $p_1$-subgroup $P_1$.  
Hence $G = P_1 \rtimes H$
for a suitable subgroup $H$ of $G$, and by Proposition 2.10 
$\psi(G)\leq \psi(P_1)\psi(H)$.
Thus
$$\psi (H) \geq \frac{\psi(G)}{\psi(P_1)}.$$
Let $|H| = h$. Then $n=h|P_1|$ and 
$\varphi (n)=\varphi (h)\varphi(|P_1|)=\varphi(h)(p_1-1)|P_1|/p_1$. Recalling 
that $p_1\geq 5$, we get
$$
\psi (H)\geq
(\frac 35)(\frac {n\varphi(n)(p_1+1)}{(p_1|P_1|^2+1)}) = (\frac 35)
(\frac {h\varphi(h)
|P_1|(p_1-1)|P_1|(p_1+1)}{p_1(p_1|P_1|^2+1)})$$
$$ > (\frac 35)(\frac {h\varphi(h) |P_1|^2(p_1^2-1)}{(p_1^2+1)|P_1|^2})=
(\frac 35)(\frac {h\varphi(h)(p_1^2-1)}{(p_1^2+1)})
\geq h\varphi(h)(\frac 35)(\frac {24}{26})
>h\varphi(h)\frac 12.$$
Let $p_2$ be the maximal prime dividing $h$. Then by Lemma 2.1 
$\psi(H) > \frac 12 \frac{h^2}{p_2}$,
and hence  there exists an element $y \in H$ satisfying
$o(y)>\frac 12\frac{h}{p_2}$. Thus
$|H: \langle y \rangle| <2p_2$ and 
by Proposition 2.5  either $H$ is solvable or there exists
a normal cyclic Sylow $p_2$-subgroup $P_2$ of $H$.  If $H$ is solvable,
then also $G$ is solvable, as required.

So suppose that there exists
a normal cyclic Sylow $p_2$-subgroup $P_2$ of $H$. Then $G = P_1
\rtimes (P_2 \rtimes V)$ for a suitable subgroup $V$ of $G$.

Now, let $p_1>p_2>\cdots > p_t>3$ be primes and suppose that
$$G = P_1 \rtimes (P_2 \rtimes (\cdots \rtimes (P_t \rtimes K))),$$
where $P_i$ are cyclic
Sylow $p_i$-subgroups of $G$ and  $K$ is a
suitable subgroup of $G$. Write $|K| = k$ and suppose that $t$ is maximal under
these conditions.
It follows from Proposition 2.10 that
$$\psi(G) \leq \psi(P_1)\psi(P_2) \cdots\psi(P_t) \psi(K), $$
and hence,  noting that $p_t>3$ and using Lemma 2.9, we get
$$\
\psi (K) \geq \frac{\psi(G)}{\psi(P_1) \psi(P_2) \cdots \psi(P_t)} \geq
\frac 35 n\varphi(n)\prod_{i=1}^t \frac {(p_i+1)}{(p_i|P_i|^2+1)}$$
$$ =
\frac 35 k\varphi(k)\prod_{i=1}^t\frac {|P_i|(p_i-1)|P_i|(p_i+1)}
{p_i(p_i|P_i|^2+1)}
=
\frac 35 k\varphi(k) \prod_{i = 1}^t
\frac{|P_i|^2(p_i^2-1)}{p_i(p_i|P_i|^2+1)}$$
$$>
\frac 35 k\varphi(k) \prod_{i = 1}^t \frac{p_i^2-1}{p_i^2+1}
>k\varphi (k)(\frac 35)(\frac 56)=
 k\varphi (k)\frac 12.$$

Let $p_{t+1}$ be the maximal prime dividing $k$. Then by Lemma 2.1 $\psi(K)>
\frac 12 \frac{k^2}{p_{t+1}},$ and there exists an element $v \in K$
satisfying
$o(v)>\frac 12\frac{k}{p_{t+1}}$. Thus
$|K: \langle v \rangle| < 2p_{t+1}$
and by Proposition 2.5 either $K$ is solvable, or there exists
a normal cyclic Sylow $p_{t+1}$-subgroup $P_{t+1}$ of $K$.  In the former case,
$K$ is solvable, and hence also $G$ is solvable, as required. If, on the 
other hand,
the latter case
occurs, then $K = P_{t+1} \rtimes  W$ for a suitable subgroup $W$ of $K$, and
by the maximality of $t$, $p_{t+1}\leq 3$. Thus $K$ is a $(2,3)$-group, hence
solvable, so also $G$ is solvable. The proof of the solvability of $G$ is
now complete.

Moreover, we have proved that $$G = P_1\rtimes (P_2 \rtimes (\dotsb (P_t
\rtimes
K))),$$
where $P_i$ are cyclic Sylow $p_i$-subgroups of $G$, and either $K$ has a cyclic maximal subgroup or 
$|K|=2^\alpha 3^\beta$ and $K$ has a cyclic subgroup of index $< 6$. We 
shall show now by induction on $t$ that these assumptions imply that
$G^{\prime\prime} \leq Z(G)$. If $t = 0$, then the
result follows from Proposition 2.4 and Proposition 2.6(2). So suppose that
$t > 0$ and set
$H =
(P_2 \rtimes \cdots (P_t \rtimes K))$.  It follows by induction 
that $H^{\prime\prime} \leq
Z(H)$. Since $G = P_1 \rtimes H$, where $P_1$ is cyclic group, 
we have $G^{\prime}
\leq
C_G(P_1)$ and $G^{\prime} = H^{\prime}[P_1, H]$. Hence 
$$G^{\prime\prime} = H^{\prime\prime}\leq Z(H) \cap C_G(P_1) \leq Z(G),$$
which completes the proof of Theorem 10.
\qed
\bigskip

Our last proof is that of Theorem 11, concerning groups of order $n$ satisfying
$\psi(G)\geq \frac 1q n\varphi(n)$.

\bigskip

\textit{Proof of Theorem 11.} Suppose that $G$ is a group of order $n$ and it 
satisfies
$\psi(G)\geq \frac 1q n\varphi(n)$.
Since by Lemma 2.1
$\varphi(n)\geq\frac {(q-1)n}p$, it follows by our assumptions that
$\psi(G)\geq \frac {(q-1)n^2}{qp}$.
Thus there exists $x\in G$
with $o(x)>\frac {(q-1)n}{qp}$ and
$$[G:\langle x\rangle]<\frac q{q-1}p\leq 2p.$$
Hence by Proposition 2.5 either $G$ has a normal cyclic
Sylow $p$-subgroup, or it is a solvable group with
a cyclic maximal subgroup of index either $p$ or $p+1$, 
as required. The proof of the theorem is now complete.
\qed

\bigskip

\section*{References}

\bibliography{mybibfile}

\end{document}